\documentclass[11pt,reqno,tbtags]{amsart}
\usepackage{amssymb}
\numberwithin{equation}{section}

\newtheorem{problem}{Problem}
\newtheorem{property}{Property}
 \let\xtheproperty\theproperty
 \def\theproperty{P\xtheproperty}
\newtheorem*{property*}{Property \csname @currentlabel\endcsname}

\newenvironment{propertyx}
{%
\begin{property*}}
{\end{property*}}

\theoremstyle{definition}

\newtheorem*{ack}{Acknowledgement}

\theoremstyle{remark}

\newenvironment{romenumerate}{\begin{enumerate}% gives (i), (ii) etc.
 }{\end{enumerate}}

\newcounter{thmenumerate}

\newcounter{xenumerate}   %no left indentation; thus wider lines

\newcommand{\refS}[1]{Section~\ref{#1}}
\newcommand{\refP}[1]{Property~\ref{#1}}

\newcommand{\refF}[1]{Figure~\ref{#1}}
\newcommand{\refand}[2]{\ref{#1} and~\ref{#2}}

\newcommand\marginal[1]{\marginpar{\raggedright\parindent=0pt\tiny #1}}

\begingroup
  \divide\count255 by 60
  \multiply\count255 by -60
  \advance\count255 by \time
  \ifnum \count255 < 10 \xdef\klockan{\the\count1.0\the\count255}
  \else\xdef\klockan{\the\count1.\the\count255}\fi
\endgroup

\newcommand\set[1]{\ensuremath{\{#1\}}}

\newcommand\Bigpar[1]{\Bigl(#1\Bigr)}

\def\rompar(#1){\textup(#1\textup)}    % usage: \rompar(...)

\newcommand\parfrac[2]{\Bigpar{\frac{#1}{#2}}}

\def\xexp(#1){e^{#1}}

\newcommand\ntoo{\ensuremath{{n\to\infty}}}

\newcommand\ie{i.e.\spacefactor=1000}
\newcommand\eg{e.g.\spacefactor=1000}

\newcommand\eqd{\overset{\mathrm{d}}{=}}

\newcounter{CC} 
\newcounter{cc}
\newcommand\E{\operatorname{\mathbb E{}}}
\renewcommand\P{\operatorname{\mathbb P{}}}
\newcommand\Var{\operatorname{Var}}

\newcommand\Po{\operatorname{Po}}
\newcommand\Bi{\operatorname{Bi}}

\newcommand\Ge{\operatorname{Ge}}

\newcommand\gss{\sigma^2}
\newcommand\eps{\varepsilon}

\newcommand\cT{{\mathcal T}}

\def\[#1]{[\![#1]\!]}

\renewcommand{\=}{:=}

\newcommand{\tn}{T_n}
\newcommand{\tni}{T_{n+1}}
\newcommand{\tx}[1]{T_{#1}}
\newcommand{\ttx}[1]{\tilde T_{#1}}
\newcommand{\too}{T_{\infty}}
\newcommand{\wk}{W_k}
\newcommand{\wx}[1]{W_{#1}}

\newcommand{\GW}{Galton--Watson}
\newcommand{\GWp}{\GW{} process}
\newcommand{\GWt}{\GW{} tree}
\newcommand{\cGWt}{conditioned \GW{} tree}

\newcommand\REM[1]{{\raggedright\texttt{[#1]}\par\marginal{XXX}}}

\newcommand\urladdrx[1]{\urladdr{\def~{\~{}}#1}}

\begin{document}
\title%[]
{Conditioned Galton--Watson trees do not grow}

\date{April 6, 2006}

\author{Svante Janson}
\address{Department of Mathematics, Uppsala University, PO Box 480,
SE-751~06 Uppsala, Sweden}
\email{svante.janson@math.uu.se}
\urladdrx{http://www.math.uu.se/~svante/}
%\urladdr{http://www.math.uu.se/\~{}svante/}

%\keywords{<keywords>}
\subjclass[2000]{} 
%{60C05 (68P10,68W40)} %%{Primary: <subject>; Secondary: <subject>}

\begin{abstract} 
An example is given which shows that, in general, \cGWt{s} cannot be
obtained by adding vertices one by one, as has been shown in a
special case by
Luczak and Winkler.
\end{abstract}

\maketitle

%\section{Introduction}\label{S:intro}
\section{Monotonicity of \cGWt{s}?}

A \cGWt{} is a random rooted tree that is (or has the same
distribution as) the family tree of a \GWp{} with some given dffspring
distribution, conditioned on the total number of vertices.

We let $\xi$ be a random variable with the given offspring
distribution; \ie, the number of offspring of each individual in the
\GWp{} is a copy of $\xi$.

We let $\xi$ be fixed throughout the paper, and let $T_n$ denote the
corresponding \cGWt{} with $n$ vertices.
For simplicity, we consider only $\xi$ such that $\P(\xi=0)>0$ and
$\P(\xi=1)>0$; then $T_n$ exists for all $n\ge1$. Furthermore, 
we assume that $\E\xi=1$ (the \GWp{} is critical) and $\gss\=\Var(\xi)<\infty$.

The importance of this construction 
lies in that 
many combinatorially interesting random trees are of this type,
for example the following:
\begin{romenumerate}
\item Random plane (= ordered) trees. $\xi\sim\Ge(1/2)$; $\gss=2$.
\item Random unordered labelled trees (Cayley trees). $\xi\sim\Po(1)$;
  $\gss=1$.
\item Random binary trees. $\xi\sim\Bi(2,1/2)$; $\gss=1/2$.
\item Random $d$-ary trees. $\xi\sim\Bi(d,1/d)$; $\gss=1-1/d$.
\end{romenumerate}
For further examples see \eg{}
Aldous \cite{AldousII} and Devroye \cite{Devroye}; 
note also that 
that the families of random trees
obtained in this way are the same as the
simply generated families of trees defined by Meir and Moon \cite{MM}.

If  we increase $n$, we get a new random tree that is in some sense
larger, but the definition above gives no relation between, say, $T_n$
and $T_{n+1}$, since they are defined by two different
conditionings. It is thus natural to ask whether 
$\tni$ is \emph{stochastically larger} than $\tn$, \ie, whether
there exists another 
construction (with the same distribution of each $T_n$) that further yields 
$T_n\subset T_{n+1}$, \ie, whether $(T_n)_{n\ge1}$ has the following
property:  
\begin{property}
  \label{P1} 
It is possible to define $T_n$ and $T_{n+1}$ on a common probability
space
such that
$T_n\subset T_{n+1}$.
\end{property}
Equivalently, \refP{P1} says that it is possible to add a new leaf to
$T_n$ by some random procedure (depending on $n$ and $\tn$) such that
the resulting tree has the distribution of $\tni$. It is thus
immediately seen that \refP{P1} is equivalent to the following:
\makeatletter\xdef\@currentlabel{\theproperty$'$}\makeatother
\begin{propertyx}
  \label{P1a}
It is possible to construct $T_1,T_2,T_3,\dots$ as a Markov chain where
at each step a new leaf is added.
\end{propertyx}

This property was investigated by Luczak and Winkler \cite{LW}, who
showed that Properties \refand{P1}{P1a} indeed hold in the case of 
random binary trees, and more generally, for random $d$-ary trees, for
any $d\ge2$.
The main purpose of this note is to give a simple counter example
(\refS{Sex}), showing that \refP{P1} does not hold for every $\xi$.

The question of whether \refP{P1} (or \ref{P1a}) holds for all \cGWt{s}
has been considered by several people, and has been explicitly stated
as an open problem at least in \cite[Problem 1.15]{SJ167}. The answer
to this question is thus negative. The problem can be reformulated as follows.

\begin{problem}
  \label{Pr1}
For which \cGWt{s} $(\tn)_n$ does \refP{P1} (or \ref{P1a}) hold?
\end{problem}

In view of the result of Luczak and Winkler \cite{LW} just mentioned, it seems
particularly interesting to study the cases of random plane trees
and random labelled trees; as far as we know, the problem is still open
for them.

It is well known that as \ntoo, $\tn$ converges in the sense of
finite-dimensional distributions to an infinite random tree $\too$
that is the family tree of the corresponding \emph{size-biased} \GWp, see
\eg{} Kennedy \cite{Kennedy},
Aldous \cite{AldousII},
Lyons, Pemantle and Peres
\cite{LPP}.
The size-biased \GWp{} 
is the same as the \emph{Q-process} studied in
\cite[Section I.14]{AN}; 
it can also
be regarded as a branching process with
two types: mortals with an offspring distribution $\xi$ and all
children mortals, and immortals with the size-biased offspring
distribution $\hat\xi$ 
with $\P(\hat\xi=j)=j\P(\xi=j)$
and exactly one immortal child (in a random
position among its siblings); the process starts with a single
immortal.
(See also \cite{SJ144}.)
Note that the 
infinite random tree $\too$  has exactly one infinite path from
the root, with (finite) \GW{} trees attached to it.

If $\xi$ is such that \refP{P1a} holds, we can construct $\tn$,
$n\ge1$, such that $T_1\subset T_2\subset\dots$, and then evidently
$T_n\to \bigcup_n\tn$; thus $\bigcup_n\tn\eqd\too$, and we may assume that
$\too=\bigcup_n\tn$. 
Hence, \refP{P1} implies the following
property: 
\begin{property}
  \label{Poo}
It is, for every $n\ge1$, possible to define $T_n$ and $\too$ on a
common probability space
such that $T_n\subset \too$.
In other words, each $\tn$ may be constructed by a suitable (random)
pruning of $\too$.
\end{property}

Thus, by Luczak and Winkler \cite{LW}, 
\refP{Poo} holds for random binary and $d$-ary trees.
On the other hand, our
counter example in \refS{Sex} also fails to satisfy \refP{Poo}.

\begin{problem}
  \label{Pr2}
For which \cGWt{s} $(\tn)_n$ does \refP{Poo} hold?
\end{problem}

Again, this problem seems to be open for random plane trees and random
labelled trees.

\section{Monotonicity of the profile?}

Properties \refand{P1}{P1a} are not only interesting in themselves,
but also technically useful (when valid), For example, for any rooted
tree $T$, let $\wk(T)$ denote the number of vertices in $T$ of
distance $k$ from the root. The sequence $(\wk(T))_{k\ge0}$ is known
as the \emph{profile} of the tree.

It is easy to see from the description of $\too$ above that 
$\E\wk(\too)=1+k\gss$. (Use the fact that
the expected number of mortal children of an immortal individual is
$\E\hat\xi-1=\E\xi^2-1=\gss$.) Moreover, as \ntoo, for each fixed $k\ge0$,
\begin{equation}
  \label{julie}
\E\wk(\tn)\to\E\wk(\too)=1+k\gss.
\end{equation}

If \refP{P1} holds, then also:
\begin{property}
  \label{Pa}
For every $k\ge0$ and $n\ge1$,
$\E \wk(\tn)\le \E\wk(\tni)$.
\end{property}
Further, if any of Property \ref{P1}, \refP{Poo} or \refP{Pa}
holds, then, using \eqref{julie}, so does the following:
\begin{property}
  \label{Pb}
For every $k\ge0$ and $n\ge1$,
$\E \wk(\tn)\le 1+k\gss$.  
\end{property}

A uniform estimate of this order, more precisely 
\begin{equation}
  \label{sofie}
\E\wk(\tn)\le Ck,
\qquad k\ge1,\,n\ge1.
\end{equation}
for all $k,n\ge1$ with a constant $C$ depending only on $\xi$, was
needed in \cite{SJ167} and proved there (Theorem 1.13) by a more
complicated argument. We will see that our counter example in
\refS{Sex} fails also \refP{Pb}; thus another argument is indeed
needed to prove \eqref{sofie} in general.

Note that Meir and Moon \cite{MM} gave explicit formulas for
$\E\wk(\tn)$ for the cases of random labelled trees, plane trees and
binary trees, which show that Properties \refand{Pa}{Pb} hold for
these cases. (Actually, the binary trees considered in \cite{MM} are
the ``strict'' or ``complete'' binary trees where all vertices have
outdegree exactly 0 or 2; these are obtained as a \cGWt{} with
$\P(\xi=0)=\P(\xi=2)=1/2$. There is a simple correspondence between
such binary trees with $2n+1$ vertices and binary trees with $n$
vertices in our notation such that, if the strict binary tree 
$\ttx{2n+1}$ corresponds to $\tn$, then 
$\wx{k+1}(\ttx{2n+1})=2\wk(\tn)$. Hence
Properties \refand{Pa}{Pb} hold for both types of random binary
trees.)

\section{A counter example}\label{Sex}

Let $\eps>0$ be a small number and let the offspring distribution be
given by
\begin{align*}
  \P(\xi=0)&=\frac{1-\eps}2,
&
  \P(\xi=1)&=\eps,
&
  \P(\xi=2)&=\frac{1-\eps}2.
\end{align*}
We have $\E\xi=1$ and $\gss\=\Var\xi=1-\eps$.
Let $\cT$ be the (unconditional) \GWt{} with this offspring
distribution. 

\begin{figure}%[htbp] 
\begin{center}
\setlength{\unitlength}{1mm} %may be changed for better appearence
\begin{picture}(20,20)
%\thicklines
\put(0,-5){\makebox(0,0){$t_1$}} 
\put(0,0){\circle*{2}} 
\put(0,10){\circle*{2}} 
\put(0,20){\circle*{2}} 
\put(0,0){\line(0,1){10}}
\put(0,10){\line(0,1){10}}
\end{picture} 
\hfil
\begin{picture}(10,10)(-10,0)
%\thicklines
\put(0,-5){\makebox(0,0){$t_2$}} 
\put(0,0){\circle*{2}} 
\put(-10,10){\circle*{2}} 
\put(10,10){\circle*{2}} 
\put(0,0){\line(1,1){10}}
\put(0,0){\line(-1,1){10}}
\end{picture} 
\end{center}
\vskip2mm
\caption{The trees with three vertices}  
\label{F1}
\end{figure}
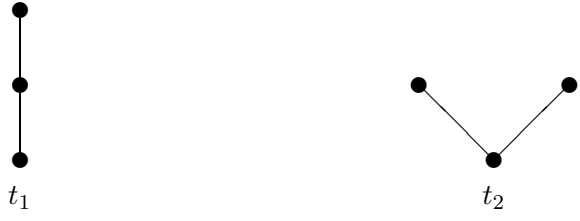

For $n=3$ we have the two possible trees in \refF{F1}.
The corresponding probabilities are, with $p_j\=\P(\xi=j)$,
\begin{align*}
  \P(\cT=t_1)&=p_1^2p_0=\eps^2\frac{1-\eps}{2}=\frac12\eps^2+O(\eps^3),
\\
  \P(\cT=t_2)&=p_2p_0^2=\parfrac{1-\eps}{2}^3=\frac18+O(\eps),
\end{align*}
and thus, conditioning on $|\cT|=3$, \ie{} on $\cT\in\set{t_1,t_2}$,
\begin{align*}
  \P(\tx3=t_1)
&
%=\P(\cT=t_1\mid |\cT|=3)
=\frac{\P(\cT=t_1)}{\P(\cT=t_1)+\P(\cT=t_2)}
=4\eps^2+O(\eps^3),
\\
  \P(\tx3=t_2)&=1-4\eps^2+O(\eps^3).
\end{align*}

\begin{figure}%[htbp] 
\begin{center}
\setlength{\unitlength}{1mm} %may be changed for better appearence
\begin{picture}(5,30)
%\thicklines
\put(0,-5){\makebox(0,0){$t_3$}} 
\put(0,0){\circle*{2}} 
\put(0,10){\circle*{2}} 
\put(0,20){\circle*{2}} 
\put(0,30){\circle*{2}} 
\put(0,0){\line(0,1){10}}
\put(0,10){\line(0,1){10}}
\put(0,20){\line(0,1){10}}
\end{picture} 
\hfil
\begin{picture}(20,20)(-10,0)
%\thicklines
\put(0,-5){\makebox(0,0){$t_4$}} 
\put(0,0){\circle*{2}} 
\put(0,10){\circle*{2}} 
\put(10,20){\circle*{2}} 
\put(-10,20){\circle*{2}} 
\put(0,0){\line(0,1){10}}
\put(0,10){\line(1,1){10}}
\put(0,10){\line(-1,1){10}}
\end{picture} 
\hfil
\begin{picture}(20,20)(-10,0)
%\thicklines
\put(0,-5){\makebox(0,0){$t_5$}} 
\put(0,0){\circle*{2}} 
\put(-10,10){\circle*{2}} 
\put(10,10){\circle*{2}} 
\put(-10,20){\circle*{2}} 
\put(0,0){\line(-1,1){10}}
\put(0,0){\line(1,1){10}}
\put(-10,10){\line(0,1){10}}
\end{picture} 
\hfil
\begin{picture}(20,20)(-10,0)
%\thicklines
\put(0,-5){\makebox(0,0){$t_6$}} 
\put(0,0){\circle*{2}} 
\put(-10,10){\circle*{2}} 
\put(10,10){\circle*{2}} 
\put(10,20){\circle*{2}} 
\put(0,0){\line(-1,1){10}}
\put(0,0){\line(1,1){10}}
\put(10,10){\line(0,1){10}}
\end{picture} 
\end{center}
\vskip3mm
\caption{The trees with four vertices}  
\label{F2}
\end{figure}
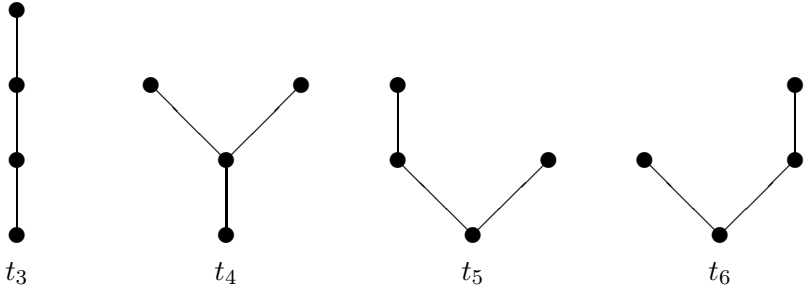

For $n=4$ we similarly have the four possible trees in \refF{F2} and
\begin{align*}
  \P(\cT=t_3)&=p_1^3p_0=\eps^3\frac{1-\eps}{2}=\frac12\eps^3+O(\eps^4),
\\
  \P(\cT=t_4)&=
  \P(\cT=t_5)=
  \P(\cT=t_6)=
p_1p_2p_0^2=\eps\parfrac{1-\eps}{2}^3=\frac18\eps+O(\eps^2),
\end{align*}
and thus, conditioning on $|\cT|=4$, %\ie{} on $\cT\in\set{t_1,t_2}$,
\begin{align*}
  \P(\tx4=t_3)
&
%=\P(\cT=t_1\mid |\cT|=3)
=O(\eps^2)
\\
  \P(\tx4=t_4)&
 = \P(\tx4=t_5)
=  \P(\tx4=t_6)
=\frac13+O(\eps^2).
\end{align*}
In particular,
\begin{align*}
  \E\wx1(\tx3)
&
=2+O(\eps^2),
\\
  \E\wx1(\tx4)
&
=\frac53+O(\eps^2),
\end{align*}
and thus $\E\wx1(\tx3)>  \E\wx1(\tx4)$ if $\eps$ is small enough, so
\refP{Pa} fails. 
(An exact calculation shows that $0<\eps<1/3$ is enough.)

By \eqref{julie}, $\E\wx1(\too)=1+\gss=2-\eps$, and thus \refP{Pb} too
fails for $k=1$, $n=3$ and small $\eps$ ($0<\eps<1/5$).
Consequently, Properties \refand{P1}{Poo} too fail.

\begin{ack}
  This research was done during a workshop at Bellairs Research
  Institute in Barbados, March 2006. I thank the participants for
  their interest and comments.
\end{ack}

\newcommand\AAP{\emph{Adv. Appl. Probab.} }
\newcommand\JAP{\emph{J. Appl. Probab.} }
\newcommand\JAMS{\emph{J. \AMS} }
\newcommand\MAMS{\emph{Memoirs \AMS} }
\newcommand\PAMS{\emph{Proc. \AMS} }
\newcommand\TAMS{\emph{Trans. \AMS} }
\newcommand\AnnMS{\emph{Ann. Math. Statist.} }
\newcommand\AnnPr{\emph{Ann. Probab.} }
\newcommand\CPC{\emph{Combin. Probab. Comput.} }
\newcommand\JMAA{\emph{J. Math. Anal. Appl.} }
\newcommand\RSA{\emph{Random Struct. Alg.} }
\newcommand\ZW{\emph{Z. Wahrsch. Verw. Gebiete} }
\newcommand\DMTCS{\jour{Discr. Math. Theor. Comput. Sci.} }

\newcommand\AMS{Amer. Math. Soc.}
\newcommand\Springer{Springer-Verlag}
\newcommand\Wiley{Wiley}

\newcommand\vol{\textbf}
\newcommand\jour{\emph}
\newcommand\book{\emph}
\newcommand\inbook{\emph}
\def\no#1#2,{\unskip#2, no. #1,} %(typeset after year) 
\newcommand\toappear{\unskip, to appear}

\newcommand\webcite[1]{\hfil\penalty0\texttt{\def~{\~{}}#1}\hfill\hfill}
\newcommand\webcitesvante{\webcite{http://www.math.uu.se/\~{}svante/papers/}}
\newcommand\arxiv[1]{\webcite{arXiv:#1}}

\def\nobibitem#1\par{}


\begin{thebibliography}{99}

%\frenchspacings

\bibitem{AldousII}
D. Aldous,
The continuum random tree II: an overview.
\emph{Stochastic Analysis (Proc., Durham, 1990)}, 23--70,
London
Math. Soc. Lecture Note Ser. 167, Cambridge Univ. Press,
Cambridge, 1991.

\bibitem{AN}
K.B. Athreya \& P.E. Ney (1972),
\book{Branching Processes}.
Springer-Verlag, Berlin, 1972.

\bibitem{Devroye}
L. Devroye,
Branching processes and their applications in the analysis of tree
structures and tree algorithms.
\inbook{Probabilistic methods for algorithmic discrete mathematics},
249--314,
eds. M. Habib et al.,
Algorithms Combin. 16, \Springer, Berlin, 1998.


\bibitem{SJ144}
S. Janson,
Ideals in a forest, one-way infinite binary trees and the contraction
method. 
\inbook{Mathematics and Computer Science II (Proceedings of the
Colloquium on Algorithms, Trees, Combinatorics and Probabilities,
Versailles 2002)}, 
393--414,
eds. B. Chauvin, P. Flajolet, D. Gardy \& A. Mokkadem, 
Trends in Mathematics, Birkh\"auser, Basel, 2002.

\bibitem{SJ167}
S. Janson,
Random cutting and records in deterministic and random trees. 
\RSA (2006),
to appear.
Available at \webcitesvante


\bibitem{Kennedy}
D.P. Kennedy (1975),
The Galton--Watson process conditioned on the total progeny.
\jour{J. Appl. Probab.} \vol{12},  800--806. 


\bibitem{LW}
M. Luczak \& P. Winkler,
Building uniformly random subtrees.
\RSA \vol{24}\no4 (2004),
420--443.

\bibitem{LPP}
R. Lyons, R. Pemantle \& Y. Peres (1995), 
Conceptual proofs of $L\log L$ criteria for mean behavior of branching
processes.
\jour{Ann. Probab.} \vol{23}, no.\ 3, 1125--1138. 

\bibitem{MM}
A. Meir \& J.W. Moon,
On the altitude of nodes in random trees.
\emph{Canad. J. Math.} \vol{30} (1978), 997--1015.

%\bibitem[??]{??}


\end{thebibliography}
\end{document}